\documentclass[10pt]{amsart}
\usepackage{enumerate,amsmath,amssymb,latexsym,
amsfonts, amsthm, amscd, mathabx}

\theoremstyle{plain}

\newtheorem{theorem}{Theorem}

\numberwithin{equation}{section}

\addtocounter{section}{-1}

\begin{document}

\title {Elliptic functions and flotation}

\date{}

\author[P.L. Robinson]{P.L. Robinson}

\address{Department of Mathematics \\ University of Florida \\ Gainesville FL 32611  USA }

\email[]{paulr@ufl.edu}

 \keywords{}

\begin{abstract}

A paraboloid or a cone of density $\rho$ with vertical axis is released from rest into a liquid of density $\rho_0$. We determine the critical value of the ratio $\rho_0/\rho$ for subsequent full submersion. 

\end{abstract}

\maketitle

\section*{Introduction} 

\medbreak 

Imagine a solid body that is able to float in a body of liquid. Now place that solid body so that its lowest point is in contact with the liquid surface and then release it; does the solid body become fully submerged? Alternatively, place the solid body so that its highest point is in contact with the liquid surface and then release it; is its buoyancy sufficient to launch it clear of the liquid? 

\medbreak 

We shall not consider these questions in full generality, but shall consider them in relation to a couple of particular examples: the case in which the solid body is a right circular cone and the case in which the solid body is cut off from a paraboloid of revolution by a plane perpendicular to its axis; as a preliminary example, we review the elementary case of a right circular cylinder. In the case of the cylinder, the vertical motion of the solid body is described by circular functions; in the other two cases, the vertical motion is described by elliptic functions. 

\medbreak 

Throughout, we treat the examples in the simplest terms, rather as problems in elementary mechanics than as realistic problems of hydrodynamics: thus, the liquid has no viscosity and such phenomena as turbulence are entirely ignored; the only forces acting vertically on the solid body are its weight and its buoyancy.

\medbreak 

\section*{Cylinder}

\medbreak 

A solid body $B$ in the form of a right circular cylinder of cross-sectional area $A$ and height $h$ is fashioned from a material of density $\rho$. This solid body $B$ is held with its axis vertical and its circular base in contact with the surface of an unlimited body of liquid having greater density $\rho_0$. We shall determine the motion of $B$ subsequent to its release from rest in this position, with the aim of deciding whether or not $B$ fully submerges. 

\medbreak 

The motion being one-dimensional, we align the downward-pointing coordinate axis along the axis of the cylinder, the origin of coordinates being the point at which this axis meets the surface of the liquid. Let $x(t)$ be the position of the centre of the circular base of $B$ at time $t$. According to the second law of Newton, if $B$ is not fully submerged then its motion is governed by the differential equation 
$$(A h \rho) x\,'' = (A h \rho) \, g - (A x \rho_0) \, g.$$
Defining  
$$\omega : = \sqrt{\frac{\rho_0 \, g}{\rho \, h}}$$
we are now presented with the initial value problem 
$$x\,'' + \omega^2 \, x = g; \; \; x(0) = 0, \, x\,'(0) = 0.$$ 
This initial value problem has solution 
$$x(t) = \frac{\rho}{\rho_0} \, h (1 - \cos \omega t)$$
where $h \, \rho/\rho_0$ is the depth at which $B$ is in equilibtium. 

\medbreak 

We are now in a position to answer the question whether or not $B$ fully submerges in the liquid; that is, whether or not $x$ achieves the value $h$ at some time. From the equation of motion, $x(t) = h$ precisely when 
$$\cos \omega t = 1 - \frac{\rho_0}{\rho} = \frac{\rho - \rho_0}{\rho}$$
which is possible precisely when $\rho_0 - \rho = |\rho - \rho_0| \leqslant \rho$ or equivalently $\rho_0 \leqslant 2 \, \rho$. We may summarize our findings as follows. If $\rho_0 > 2 \, \rho$ then $B$ does not fully submerge; the base of $B$ reaches a maximum depth $2 \, h \, \rho / \rho_0$ at time $\pi / \omega$. If $\rho_0 < 2 \, \rho$ then $B$ becomes fully submerged, at which point in time the differential equation governing the motion changes: where $x > h$ the buoyancy is constant and $x\,'' = (1 - \frac{\rho_0}{\rho} )\, g$. In the borderline case $\rho_0 = 2 \rho$ the cylinder $B$ falls until its circular top is level with the liquid surface and then reverses direction. 

\medbreak 

Suppose instead that the cylinder $B$ is again initially at rest but has its circular top flush with the liquid surface. The motion of $B$ is again governed by the differential equation $x\,'' + \omega^2 \, x = g$ but now the initial conditions are $x(0) = h$ and $x \, '(0) = 0$. This new initial value problem has solution 
$$x(t) = \frac{\rho}{\rho_0} \, h + (1 - \frac{\rho}{\rho_0}) \, h \cos \omega t.$$ 
We may now answer the question whether or not the upward motion of $B$ launches it clear of the liquid; that is, whether or not $x$ achieves the value $0$ at some time. Again, the critical value of the ratio $\rho_0$ to $\rho$ is $2$: if this ratio is less than $2$ then the cylinder remains partly submerged; if it is greater than $2$ then the cylinder leaves the liquid (and is then subject to gravity alone). 

\medbreak 

\section*{Paraboloid}

\medbreak 

We now turn to the more interesting case in which the solid body $B$ is cut from a paraboloid of revolution. To be specific, in cylindrical polar coordinates, let the solid paraboloid be given by `$r^2 \leqslant 2 p z$' and let $B$ be the bounded portion that is cut off by the plane `$z = h$'. The volume of $B$ is then 
$$\int_0^h \pi r^2 {\rm d} z = \int_0^h \pi 2 p z {\rm d} z = \pi p h^2.$$

\medbreak 

Orient $B$ so that its axis is vertical, its circular surface being above its vertex; place the vertex of the solid on the surface of the liquid and release it from rest. As before, we wish to know whether or not $B$ becomes fully submerged during the course of its descent into the liquid, assuming that the solid has density $\rho$ and the liquid has density $\rho_0$.  

\medbreak 

Once again, we let the axis of $B$ be the downward-pointing coordinate axis, taking the point at which it meets the surface of the liquid to be the origin. Let $x(t)$ be the position of the vertex of $B$ at time $t$: thus, $x$ is the depth at which the vertex of the partly-submerged $B$ lies below the surface of the liquid. The mass of $B$ is $\pi p h^2 \rho$; when the vertex of $B$ is at depth $x$ below the surface, mass $\pi p x^2 \rho_0$ of liquid is displaced. So long as $B$ is not fully submerged, its motion is therefore governed by the nonlinear differential equation 
$$(\pi p h^2 \rho) x\,'' = (\pi p h^2 \rho) g - (\pi p x^2 \rho_0) g.$$
We are thus presented with the initial value problem 
$$x\,'' = g - \frac{\rho_0 \, g}{\rho \, h^2} x^2; \; \; x(0) = 0, \, x\,'(0) = 0$$
which has first integral
$$\frac{1}{2} (x\,')^2 = g \, x - \frac{\rho_0 \, g}{\rho \, h^2} \frac{x^3}{3}; \; x(0) = 0.$$

\medbreak 

We are now prepared to decide whether or not $B$ becomes fully submerged; this decision can be reached by direct inspection of the first integral. As $B$ falls, $x$ becomes more positive until $x \, ' = 0$ and 
$$x =  h \sqrt{\tfrac{3 \rho}{\rho_0}}$$
unless full submersion has already taken place. Consequently, the critical value of the ratio $\rho_0 / \rho$ is now $3$. If $\rho_0 > 3 \rho$ then the maximum value of $x$ is less than $h$ and $B$ does not fully submerge. If $\rho_0 < 3 \rho$ then full submersion of $B$ takes place, whereupon the buoyancy levels off and the differential equation for $x$ changes. 

\medbreak 

We wish to recast the differential equation and its first integral into more familiar terms. For this purpose, we shall reverse the sign of $x$ and incorporate appropriate dimensional parameters. Explicitly, introduce a rescaling $u$ of $x$ by the rule 
$$x(t) = - \lambda \, u\,(\omega  \, t)$$ 
where $\lambda$ and $\omega$ are positive constants (with dimensions $L$ and $1/T$ respectively). The differential equation satisfied by $x$ is thereby converted to the differential equation 
$$u\,'' = \frac{\rho_0 \, g}{\rho \, h^2} \frac{\lambda}{\omega^2} \, u^2 - \frac{g}{\lambda \, \omega^2}$$
for $u$. Let us specify the values of the dimensional parameters by 
$$\lambda^2 = 12 \, \frac{\rho}{\rho_0} h^2$$ 
and 
$$\omega^4 = \frac{1}{3} \frac{\rho_0}{\rho} \frac{g^2}{h^2}.$$
The differential equation for $u$ then becomes the familiar 
$$u\,'' = 6 u^2 - \tfrac{1}{2}$$
with first integral 
$$(u\,')^2 = 4 u^3 - u - g_3$$ 
where $g_3$ is a constant, which the initial conditions $u(0) = 0$ and $u\,'(0) = 0$ force to be $0$. Conclusion: the rescaled position $u$ is the (lemniscatic) Weierstrass $\wp$-function, with invariants $g_2 = 1$ and $g_3 = 0$, translated so that the origin is not a pole but rather a zero. 

\medbreak 

Incidentally, we may answer afresh the question whether or not $B$ becomes fully submerged from knowledge of the lemniscatic $\wp$-function. Let $a$ be a zero of this $\wp$-function; as $z$ varies along the line ${\rm Im} \, z = {\rm Im} \, a$, the value of the $\wp$-function at $z$ oscillates between $0$ and $- 1/2$. Note that  value $-1/2$ of the shifted $\wp$-function $u$ corresponds to the value $h \sqrt{3 \rho/\rho_0}$ of $x$. 

\medbreak 

The question whether or not $B$ is launched clear of the liquid, if it is initially at rest with its circular top at the surface of the liquid, is left as an exercise. We merely remark that $u$ is again a shifted $\wp$-function for which $g_2 = 1$ but for which 
$$g_3 = 4 u(0)^3 - u(0) - u\,'(0)^2 = \frac{1}{2} \sqrt{\frac{\rho_0}{3 \rho}}\, \Big(1 - \frac{\rho_0}{3 \rho}\Big).$$

\medbreak 

\section*{Cone}

\medbreak 

The case in which the solid body $B$ is a right circular cone also involves elliptic functions. 

\medbreak 

Let $B$ have height $h$ and circular area $A$; again let the solid have density $\rho$ and consider its flotation in a liquid of density $\rho_0$. Initially, suspend the cone with its axis vertical and with only its vertex in contact with the liquid; then release the cone from rest. We again wish to determine whether or not the solid $B$ fully submerges in the liquid as it descends. 

\medbreak 

Let $x$ denote the depth of the vertex of $B$ below the surface of the liquid. Provided that $B$ is only partly submerged, its motion follows the differential equation 
$$\tfrac{1}{3} A h \rho \, x\,'' = \tfrac{1}{3} A h \rho \, g - \tfrac{1}{3} \big(A \frac{x^2}{h^2}\big) x \rho_0 \, g$$
or 
$$x\,'' = g - \frac{\rho_0 \, x^3}{\rho \, h^3} g$$
with initial conditions $x(0) = 0$ and $x\,'(0) = 0$. This initial value problem has first integral 
$$\tfrac{1}{2} (x\, ')^2 = g \, x - \tfrac{1}{4} \frac{\rho_0 \, g}{\rho \, h^3} x^4$$
or 
$$(x\,')^2 = \frac{\rho_0 \, g}{2 \rho \, h^3} \, x \, \big(\frac{4 \rho}{\rho_0} h^3 -x^3 \big)$$
with initial condition $x(0) = 0$. 

\medbreak 

As before, the question whether or not $B$ becomes fully submerged during its fall may be answered directly from the differential equation: $B$ falls until $x\,' = 0$ when 
$$x = (\frac{4 \rho}{\rho_0})^{1/3} h$$ 
unless full submersion has already taken place; the crucial value of the ratio $\rho_0/\rho$ is thus $4$. If $\rho_0 > 4 \rho$ then the maximum value of $x$ displayed above is less than $h$ and $B$ does not fully submerge. If $\rho_0 < 4 \rho$ then $B$ fully submerges before this maximum can be reached. 

\medbreak

Unless and until full submersion occurs, $x$ is again an elliptic function. To identify this elliptic function in familiar terms, we rescale. Introduce new dimensional parameters $\lambda > 0$ and $\omega > 0$ by 
$$\lambda = (\frac{4 \rho}{\rho_0})^{1/3} h$$ 
and 
$$\omega^2 = \frac{2 g}{\lambda} = \frac{2 g}{h} (\frac{\rho_0}{4 \rho})^{1/3} .$$ 
The rescaled function $u$ given by 
$$x(t) = \lambda \, u(\omega \, t)$$ 
satisfies the first-order differential equation 
$$(u\,')^2 = u \,(1 - u^3).$$ 
The reciprocal $v = 1/u$ then satisfies 
$$(v\,')^2 = v^3 - 1$$
and the further rescaling $w = v/4$ leads to 
$$(w\,')^2 = 4 w^3 - \frac{1}{16}.$$ 
As the position $x$ has a zero at the origin, it follows that the function $w$ has a pole at the origin. Conclusion: $w$ is the Weierstrass $\wp$-function with invariants $g_2 = 0$ and $g_3 = 1/16$; further rescaling will make $g_3 = 1$.  

\medbreak 

\section*{Remarks}

\medbreak 

More generally, we may take the solid body $B$ to be the bounded region cut off from a solid of revolution by a plane perpendicular to its axis. For example, suppose that $B$ is cut off from the region $0 \leqslant r \leqslant k z^d$ by the plane $z = h$ (in cylindrical polar coordinates). As before, suppose that $B$ has density $\rho$ and consider its immersion in a body of liquid having density $\rho_0$. Initially, orient $B$ so that its axis is vertical and only its vertex is in contact with the liquid; then release it from rest. Unless and until full submersion occurs, the depth $x$ of the vertex of $B$ below the surface of the liquid satisfies
$$\frac{x\,''}{g} = 1 - \frac{\rho_0}{\rho} \frac{x^{2 d + 1}}{h^{2 d + 1}}; \; \; x(0) = 0, \, x\,'(0) = 0$$
with first integral 
$$\frac{(x\,')^2}{2 g} = x \, \Big(1 - \frac{\rho_0}{2 (d + 1) \rho} \frac{x^{2 d + 1}}{h^{2 d + 1}} \Big); \; x(0) = 0.$$
The critical value of the ratio $\rho_0/\rho$ is therefore $2(d + 1)$ for full submersion. Here, the cases in which $2 d$ is $0, 1$ and $ 2$ correspond respectively to the cylinder, the paraboloid and the cone. Greater values of $d$ lead to hyperelliptic functions prior to full submersion. 

\medbreak 

We may also approach certain other solids by similar methods. The case of a right equilateral triangular prism stands out: if such a prism (for which $h$ is the height of a triangular section) is released from rest when one lateral face is horizontal and the opposite lateral edge is its only contact with the liquid, then its equation of motion is precisely the same as that previously found for the paraboloid of revolution; in particular, the critical density ratio for full submersion is $\rho_0/\rho = 3$. If instead the axis of the prism is oriented vertically, then the prism behaves as the cylinder with which we began our account, so the critical density ratio is $2$. 

\medbreak 

The appearance of elliptic functions in the motion of a solid body through a liquid medium is discussed in the Appendix to the classic Greenhill [1] and references therein: there, the motion of a fully submerged solid of revolution is considered quite generally; here, we focus more simply the process of submersion. The more modern Lawden [2] presents further appearances of elliptic functions in physical problems, in Chapter 5 and Chapter 7, though none of these appearances is directly related to hydrodynamics. 

\medbreak

\bigbreak 

\begin{center} 
{\small R}{\footnotesize EFERENCES}
\end{center} 
\medbreak 

[1] G. Greenhill, {\it The Applications of Elliptic Functions}, Macmillan and Company (1892). 

\medbreak 

[2] D. F. Lawden, {\it Elliptic Functions and Applications}, Springer-Verlag (1989). 

\medbreak

\end{document}